# Some beautiful q- analogues of Fibonacci and Lucas polynomials

*Johann Cigler*

Fakultät für Mathematik, Universität Wien

johann.cigler@univie.ac.at

## Abstract

We give an overview about well-known basic properties of two classes of $q-$Fibonacci and $q-$Lucas polynomials and offer a common generalization.

## 1. Introduction

We study two classes of $q-$Fibonacci polynomials with interesting properties: the Carlitz $q-$Fibonacci polynomials $F_n(x,s,q)$ and the polynomials $Fib_n(x,s,q)$ introduced in [6] and [8].

More generally we consider the polynomials

$$\Phi_n(x,s,m,q)=\sum_{k=0}^{\left\lfloor\frac{n-1}{2}\right\rfloor} q^{\binom{k+1}{2}+m\binom{k}{2}} \begin{bmatrix} n-1-k \\ k \end{bmatrix} s^k x^{n-1-2k} \tag{1.1}$$

for $m\in\mathbb{Z}$ which reduce to $F_n(x,s,q)$ for $m=1$ and to $Fib_n(x,s,q)$ for $m=0$.

They satisfy the recurrence

$$\Phi_n(x,s,m,q)=x\Phi_{n-1}(x,qs,m,q)+qs\Phi_{n-2}(x,q^{m+1}s,m,q) \tag{1.2}$$

with initial values $\Phi_0(x,s,m,q)=0$ and $\Phi_1(x,s,m,q)=1$.

The associated $q-$Lucas polynomials $\Lambda_n(x,s,m,q)$ are given by

$$\Lambda_n(x,s,m,q)=\sum_{k=0}^{\frac{n}{2}} q^{(m+1)\binom{k}{2}} \frac{[n]}{[n-k]} \begin{bmatrix} n-k \\ k \end{bmatrix} s^k x^{n-2k}. \tag{1.3}$$

Although these do not satisfy simple recurrences such as (1.2) there is another type of recurrence which we will call $D-$recurrence which is satisfied by both $q-$Fibonacci and $q-$Lucas polynomials: Let $D$ denotes the $q-$differentiation operator. Then

$$\Phi_n(x,s,m,q) = x\Phi_{n-1}(x,s,m,q) + (q-1)sD\Phi_{n-1}(x,q^m s,m,q) + s\Phi_{n-2}(x,q^m s,m,q) \quad (1.4)$$

and

$$\Lambda_n(x,s,m,q) = x\Lambda_{n-1}(x,s,m,q) + (q-1)sD\Lambda_{n-1}(x,q^m s,m,q) + s\Lambda_{n-2}(x,q^m s,m,q). \quad (1.5)$$

Of course for $q=1$ this reduces to the recursions $F_n(x,s) = xF_{n-1}(x,s) + sF_{n-2}(x,s)$ and $L_n(x,s) = xL_{n-1}(x,s) + sL_{n-2}(x,s)$ of the classical Fibonacci and Lucas polynomials.

First of all we sketch the roots of the theory which are well-known properties of Fibonacci and Lucas numbers and polynomials. Then we recall the main results about the two different $q-$ analogues and finally we expose some facts which give rise to the above mentioned common generalization.

## 2. The concrete roots of the theory

### 2.1.

Let us begin with the **Fibonacci numbers** $F_n$. The first few terms are $0,1,1,2,3,5,8,13,21,34,55,\cdots$.

They satisfy

$$F_n = F_{n-1} + F_{n-2} \quad (2.1)$$

for $n \geq 2$ with initial values $F_0 = 0$ and $F_1 = 1$.

#### 2.1.1.

There are many methods how to deal with Fibonacci numbers. The simplest one is to find numbers $x$ which satisfy $x^n = x^{n-1} + x^{n-2}$ and write the Fibonacci numbers as linear combinations of these numbers. This leads to **Binet's formula**

$$F_n = \frac{\alpha^n - \beta^n}{\alpha - \beta}, \quad (2.2)$$

where $\alpha = \frac{1+\sqrt{5}}{2}$ and $\beta = \frac{1-\sqrt{5}}{2}$ are the solutions of $x^2 - x - 1 = 0$.

Formula (2.2) can be used to extend $F_n$ to negative indices and yields

$$F_{-n} = (-1)^{n-1} F_n. \tag{2.3}$$

**2.1.2.**

Another method looks for powers of matrices with the same recurrence relation:

Let $C = \begin{pmatrix} 0 & 1 \\ 1 & 1 \end{pmatrix}$. Then $C^2 = \begin{pmatrix} 0 & 1 \\ 1 & 1 \end{pmatrix}^2 = \begin{pmatrix} 1 & 1 \\ 1 & 2 \end{pmatrix} = \begin{pmatrix} 1 & 0 \\ 0 & 1 \end{pmatrix} + \begin{pmatrix} 0 & 1 \\ 1 & 1 \end{pmatrix} = I + C,$

where $I$ denotes the identity matrix.

The matrix powers are easily computed to be

$$C^n = \begin{pmatrix} F_{n-1} & F_n \\ F_n & F_{n+1} \end{pmatrix}. \tag{2.4}$$

The traces of these matrices are the **Lucas numbers** $L_n$,

$$L_n = tr\left(C^n\right) = F_{n-1} + F_{n+1} = \alpha^n + \beta^n. \tag{2.5}$$

They satisfy the same recurrence as the Fibonacci numbers, but with initial values

$L_0 = 2$ and $L_1 = 1$. The first terms are $2,1,3,4,7,11,18,29,47,\cdots$.

The determinants of these matrices give **Cassini's formula**

$$F_{n-1}F_{n+1} - F_n^2 = (-1)^n. \tag{2.6}$$

From $C^{m+n} = C^m C^n$ we deduce

$$\begin{pmatrix} F_{m+n-1} & F_{m+n} \\ F_{m+n} & F_{m+n+1} \end{pmatrix} = \begin{pmatrix} F_{m-1} & F_m \\ F_m & F_{m+1} \end{pmatrix} \begin{pmatrix} F_{n-1} & F_n \\ F_n & F_{n+1} \end{pmatrix}$$

and identities such as

$$F_{m+n} = F_m F_{n-1} + F_{m+1} F_n, \tag{2.7}$$

with interesting arithmetical consequences.

For $n = m$ this formula connects Fibonacci and Lucas numbers

$$F_{2n} = F_n F_{n-1} + F_{n+1} F_n = F_n (F_{n-1} + F_{n+1}) = F_n L_n. \tag{2.8}$$

**2.1.3.**

There are many **combinatorial interpretations** of the Fibonacci numbers. We choose the following one: Consider subintervals of the integers of the form $\{m, m+1, \cdots, m+p-1\}$ and cover these sets with dots $a$ and dashes $b.$ A dot $a$ covers one number and a dash $b$ covers two adjacent numbers. A covering $c = c_1 c_2 \cdots c_r$ of $\{m, m+1, \cdots, m+p-1\}$ with $c_i \in \{a, b\}$ is called a **Fibonacci word** or Morse sequence **of length** $p.$ For example let $m = 2$ and $p = 5.$ Then there are 8 Fibonacci words $aaaaa, aaab, aaba, abaa, abb, baaa, bab, bba$ on $\{2,3,4,5,6\}$.

The number of Fibonacci words of length $n-1$ is equal to $F_n.$ This is obvious for $n \geq 2$ since there is a single Fibonacci word $a$ on $\{0\}$ and there are $2$ Fibonacci words on $\{0,1\},$ namely $aa$ and $b.$ In the general case a Fibonacci word of length $n-1$ either begins with a dot with $F_{n-1}$ possibilities for the remaining word or it begins with a dash, in which case there are $F_{n-2}$ possibilities for the remaining word. The cases $n = 0$ and $n = 1$ can be subsumed under this interpretation by pretending that for $n = 1$ there is exactly one Fibonacci word, the empty word $\varepsilon$ and that for $n = 0$ there exists no Fibonacci word.

The number $f(n,k)$ of Fibonacci words of length $n-1$ with $k$ dashes is $f(n,k) = \binom{n-1-k}{k}.$ For in this case the word consists of $k + (n-1-2k) = n-1-k$ elements and has $\binom{n-1-k}{k}$ possibilities for the $k$ dashes.

**2.1.4.**

We will also need an **algebraic version** of this interpretation. To this end consider all finite linear combinations with complex coefficients of the words
$\varepsilon, a, b, aa, ab, ba, bb, aaa, aab, aba, abb, baa, bab, bba, bbb, \cdots.$

If we define the product of two words $c = c_1 \cdots c_m$ and $d = d_1 \cdots d_n$ as the juxtaposition $c_1 \cdots c_m d_1 \cdots d_n$ we get the ring $P(a,b)$ of all polynomials in two non-commuting letters $a, b.$

Let $C_k^n(a,b)$ be the sum of all words with $k$ letters $b$ and $n-k$ letters $a.$

Then obviously

$$(a+b)^n = \sum_{k=0}^{n} C_k^n(a,b) \tag{2.9}$$

For example

$$(a+b)^2 = (a+b)(a+b) = a^2 + ba + ab + b^2 = C_0^2(a,b) + C_1^2(a,b) + C_2^2(a,b).$$

Define $F_n^*(a,b) \in P(a,b)$ as the sum of all Fibonacci words of length $n-1$. The sequence $\left(F_n^*(a,b)\right)$ begins with $0, \varepsilon, a, a^2+b, a^3+ab+ba, a^4+a^2b+aba+ba^2+b^2, \cdots$.

These polynomials satisfy

$$F_n^*(a,b) = aF_{n-1}^*(a,b) + bF_{n-2}^*(a,b) \tag{2.10}$$

and also

$$F_n^*(a,b) = F_{n-1}^*(a,b)a + F_{n-2}^*(a,b)b. \tag{2.11}$$

It is clear that

$$F_n^*(a,b) = \sum_{k=0}^{n-1} C_k^{n-1-k}(a,b). \tag{2.12}$$

For example

$$F_5^*(a,b) = C_0^4(a,b) + C_1^3(a,b) + C_2^2(a,b) = aaaa + \left(aab + aba + baa\right) + \left(bb\right) = a^4 + a^2b + aba + ba^2 + b^2.$$

**2.1.5.**

The Lucas numbers can be interpreted by using "circular" Fibonacci words.

Divide the circumference of a circle into $n$ arcs with equal lengths, denote them by $0, 1, \cdots, n-1$ and cover them with dots and dashes. Let $l(n,k)$ be the number of coverings with $k$ dashes. Then $l(n,k) = f(n+1,k) + f(n-1,k-1) = \binom{n-k}{k} + \binom{n-1-k}{k-1} = \frac{n}{n-k}\binom{n-k}{k}$.

This identity comes from associating to a circular word a linear word by a suitable cut.

If the arc $0$ is a dot or the first point of a dash we cut the word before this place and get a linear word of length $n$ with $k$ dashes. If it is the end of a dash we eliminate the dash and there remains a linear word of length $n-2$ with $k-1$ dashes. Thus $L_n = F_{n+1} + F_{n-1}$ in accord with (2.5).

For $n=0$ we must set $l(0,0)=2$ in order to have $L_0=2.$ I did not look for a logical sophism to "explain" this choice.

### 2.2.

The combinatorial interpretation of the Fibonacci numbers leads at once to the **Fibonacci polynomials**: We associate with each Fibonacci word consisting of $k$ dashes and $n-1-2k$ dots the weight $s^k x^{n-1-2k}$. Since there are $f(n,k)=\binom{n-1-k}{k}$ such words the weight of all Fibonacci words of length $n-1$ is

$$F_n(x,s)=\sum_{k=0}^{\left\lfloor\frac{n-1}{2}\right\rfloor}\binom{n-1-k}{k}s^k x^{n-1-2k}. \tag{2.13}$$

Classify all Fibonacci words with respect to their first letter. The weight of all words with first letter $a$ is $xF_{n-1}(x,s)$ and the weight of those beginning with $b$ is $sF_{n-2}(x,s)$ ). This gives the recurrence

$$F_n(x,s)=xF_{n-1}(x,s)+sF_{n-2}(x,s) \tag{2.14}$$

*for* $n\geq 2$ *with initial values* $F_0(x,s)=0$ *and* $F_1(x,s)=1.$

#### 2.2.1.

Translated into algebraic language this can be formulated in the following way: Consider the ring homomorphism $\varphi$ from $P(a,b)$ into the polynomials defined by $\varphi(a)=x$ and $\varphi(b)=s.$ Then $\varphi\left((a+b)^n\right)=(x+s)^n,\ \ \varphi\left(C_k^n(a,b)\right)=\binom{n}{k}s^k x^{n-k}$ and

$$\varphi\left(F_n^*(a,b)\right)=\varphi\left(\sum_{k=0}^{\left\lfloor\frac{n-1}{2}\right\rfloor}C_k^{n-1-k}(a,b)\right)=\sum_{k=0}^{\left\lfloor\frac{n-1}{2}\right\rfloor}\binom{n-1-k}{k}s^k x^{n-1-2k}=F_n(x,s).$$

**2.2.2.**

Binet's formula becomes

$$F_n(x,s)=\frac{\alpha^n-\beta^n}{\alpha-\beta}, \tag{2.15}$$

where $\alpha=\frac{x+\sqrt{x^2+4s}}{2}$ and $\beta=\frac{x-\sqrt{x^2+4s}}{2}$ are the solutions of $z^2-xz-s=0.$

This formula can be used to extend $F_n(x,s)$ to negative indices and yields

$$F_{-n}(x,s)=(-1)^{n-1}\frac{F_n(x,s)}{s^n}. \tag{2.16}$$

To motivate a later generalization let us note that $\left(F_n(1,-1)\right)_{n\geq 0}=\left(1,0,-1,-1,0,1,\cdots\right)$ with period $6.$ The reason is that $\alpha=\frac{1+\sqrt{-3}}{2},\beta=\frac{1-\sqrt{-3}}{2}$ are sixth roots of unity.

**2.2.3.**

The matrix $C(x,s)=\begin{pmatrix}0 & 1\\ x & s\end{pmatrix}$ satisfies

$$C(x,s)^2=\begin{pmatrix}0 & 1\\ x & s\end{pmatrix}^2=\begin{pmatrix}x & s\\ xs & x+s^2\end{pmatrix}=x\begin{pmatrix}1 & 0\\ 0 & 1\end{pmatrix}+s\begin{pmatrix}0 & 1\\ x & s\end{pmatrix}=xI+sC(x,s).$$

The powers give

$$C(x,s)^n=\begin{pmatrix}sF_{n-1}(x,s) & F_n(x,s)\\ sF_n(x,s) & F_{n+1}(x,s)\end{pmatrix}. \tag{2.17}$$

The traces of these matrices are the **Lucas polynomials** $L_n(x,s)$

$$L_n(x,s)=tr\left(C(x,s)^n\right)=sF_{n-1}(x,s)+F_{n+1}(x,s)=\alpha^n+\beta^n. \tag{2.18}$$

They satisfy the same recurrence as the Fibonacci polynomials, but with initial values

$L_0(x,s)=2$ and $L_1(x,s)=x.$

An explicit formula for $n>0$ is

$$L_n(x,s)=\sum_{k=0}^{\left\lfloor \frac{n}{2}\right\rfloor}\frac{n}{n-k}\binom{n-k}{k}s^k x^{n-2k}. \qquad (2.19)$$

From $C(x,s)^{m+n}=C(x,s)^m C(x,s)^n$ or

$$\begin{pmatrix} sF_{m-1}(x,s) & F_m(x,s) \\ sF_m(x,s) & F_{m+1}(x,s)\end{pmatrix}\begin{pmatrix} sF_{n-1}(x,s) & F_n(x,s) \\ sF_n(x,s) & F_{n+1}(x,s)\end{pmatrix}=\begin{pmatrix} sF_{m+n-1}(x,s) & F_{m+n}(x,s) \\ sF_{m+n}(x,s) & F_{m+n+1}(x,s)\end{pmatrix}$$

we conclude

$$F_{m+n}(x,s)=F_m(x,s)F_{n+1}(x,s)+sF_n(x,s)F_{m-1}(x,s). \qquad (2.20)$$

For $m=n$ this gives

$$F_{2n}(x,s)=F_n(x,s)\left(F_{n+1}(x,s)+sF_{n-1}(x,s)\right)=F_n(x,s)L_n(x,s). \qquad (2.21)$$

The determinants of these matrices give Cassini's formula for Fibonacci polynomials

$$F_{n-1}(x,s)F_{n+1}(x,s)-\left(F_n(x,s)\right)^2=(-1)^n s^{n-1}. \qquad (2.22)$$

**2.2.4.**

The Lucas polynomials can again be interpreted by using circular Fibonacci words as above.

The weight of all words where the arc $0$ is the endpoint of a dash is $sF_{n-1}(x,s)$ and the weight of the other words is $F_{n+1}(x,s)$.

This also gives a combinatorial interpretation of the formula

$$L_n(x,s)=F_{n+1}(x,s)+sF_{n-1}(x,s). \qquad (2.23)$$

For some purposes it is useful to consider **modified Lucas polynomials** $L_n^*(x,s)$ which coincide with $L_n(x,s)$ for $n>0$ but with initial value $L_0^*(x,s)=1$.

They satisfy

$$\sum_{k=0}^{\left\lfloor \frac{n}{2}\right\rfloor}(-s)^k\binom{n}{k}L_{n-2k}^*(x,s)=x^n. \qquad (2.24)$$

We prove this by induction. It is obviously true for $n=0$ and $n=1$. If we set $L_n^*(x,s)=0$ for $n<0$ then the assertion follows from the following computation

$$\sum_k (-s)^k \binom{n+1}{k} L_{n+1-2k}^*(x,s) = \sum_k (-s)^k \left( \binom{n}{k} + \binom{n}{k-1} \right) L_{n+1-2k}^*(x,s)$$
$$= \sum_k (-s)^k \binom{n}{k} L_{n+1-2k}^*(x,s) + \sum_k (-s)^{k+1} \binom{n}{k} L_{n-1-2k}^*(x,s)$$
$$= \sum_k (-s)^k \binom{n}{k} \left( L_{n+1-2k}^*(x,s) - sL_{n-1-2k}^*(x,s) \right) = x \sum_k (-s)^k \binom{n}{k} L_{n-2k}^*(x,s).$$

We have only to observe that $L_{k+1}^*(x,s) - sL_{k-1}^*(x,s) = xL_k^*(x,s)$ for $k \neq 1$ and that $L_2^*(x,s) - sL_0^*(x,s)$ can only occur if $n = 2m+1$ is odd and $k = m$. In this case we get

$$(-s)^m \binom{2m+1}{m} \left( L_2^*(x,s) - sL_0^*(x,s) \right) + (-s)^{m+1} \binom{2m+1}{m+1} = (-s)^m \binom{2m+1}{m} x^2 = xL_1^*(x,s).$$

## Remarks

There is some confusion in the literature about the numbering of Fibonacci numbers and Fibonacci polynomials. We would have avoided some logical troubles if we would have considered $f_n = F_{n+1}$ instead of $F_n$ as many authors do. A further argument for this choice would be that the degree of $f_n(x,s) = F_{n+1}(x,s)$ as polynomial in $x$ would be the same as the index. But there is an important number theoretical property of Fibonacci numbers and polynomials which is only true for $F_n$. This is the fact that $F_{kn}$ is a multiple of $F_n$ and more generally

$$\gcd\left(F_m, F_n\right) = F_{\gcd(m,n)}. \tag{2.25}$$

Therefore I prefer the above definition because it does not hide these underlying mathematical facts.

It should be noted that there are also close connections with **Chebyshev polynomials**. These are defined by $T_n(\cos\vartheta) = \cos n\vartheta$ and $U_n(\cos\vartheta) = \dfrac{\sin(n+1)\vartheta}{\sin\vartheta}$ and satisfy $T_n(x) = \dfrac{L_n(2x,-1)}{2}$ and $U_n(x) = F_{n+1}(2x,-1).$

Some identities for Fibonacci and Lucas polynomials can be interpreted as generalizations of **trigonometric identities**. For example the identity $L_n(x,s)^2 - (x^2+4s)F_n(x,s)^2 = 4(-s)^n$ generalizes the trigonometric identity $\cos^2 n\vartheta + \sin^2 n\vartheta = 1.$ From (2.15) and (2.18) follows that

$$\alpha^n = \frac{L_n(x,s) + F_n(x,s)\sqrt{x^2+4s}}{2} = \frac{L_n(x,s) + \sqrt{(x^2+4s)F_n(x,s)^2}}{2} = \frac{L_n(x,s) + \sqrt{L_n(x,s)^2 - 4(-s)^n}}{2}$$

which implies $L_{mn}(x,-s) = L_m(L_n(x,-s),-s^n).$ This reduces to the **composition identity** $T_{mn}(x) = T_m\left(T_n(x)\right)$ for the Chebyshev polynomials of the first kind if we let $x \to 2x$ and $s=1.$ It would be very interesting to find $q-$analogues of these formulas.

There are two simple ways to use our combinatorial model to obtain $q-$analogues of the Fibonacci polynomials.

1) For a given Fibonacci word $c = c_1c_2\cdots c_{n-1-k}$ on $\{0,1,\cdots,n-2\}$ let $1 \le j_1 < j_2 < \cdots < j_k \le n-2$ be the endpoints of the $k$ dashes. Associate with $c$ the weight $w(c) = q^{j_1+j_2+\cdots+j_k}s^kx^{n-1-2k}$ and define $F_n(x,s,q) = \sum_c w(c)$ where the sum extends over all Fibonacci words of length $n-1.$

2) Let $1 \le i_1 < i_2 < \cdots < i_k \le n-1-k$ be the indices such that $c_{i_j} = b$. Associate with $c$ the weight $W(c) = q^{i_1+i_2+\cdots+i_k}s^kx^{n-1-2k}$ and define $Fib_n(x,s,q) = \sum_c W(c).$

It turns out that both choices are "natural" for some applications.

## 3. Carlitz's q-Fibonacci and q-Lucas polynomials

### 3.1.

To each Fibonacci word $c$ on $\{m, m+1,\cdots, m+p-1\}$ we associate a **weight** $w(c)$ in the following way: If the number $i$ is covered by a dot then $w(i) = x,$ if it is the endpoint of a dash we set $w(i) = q^is$ and $w(i) = 1$ in all other cases. Then $w(c) = \prod_{i=m}^{m+p-1} w(i).$

Thus the weights of the 8 Fibonacci words $aaaaa, aaab, aaba, abaa, abb, baaa, bab, bba$ on $\{2,3,4,5,6\}$ are $x^5, x^3q^6s, x^3q^5s, x^3q^4s, xq^{10}s^2, x^3q^3s, xq^9s^2, xq^8s^2.$

**Definition 3.1**

*The* ***(Carlitz-)*** $q-$***Fibonacci polynomials*** $F_n(x,s,q)$ *are defined as*

$$F_n(x,s,q)=\sum_c w(c) \tag{3.1}$$

*where the sum is extended over all Fibonacci words* $c$ *on* $\{0,1,\cdots,n-2\}.$

The $q-$Fibonacci polynomials $F_n(x,s,q)$ satisfy the recursion

$$F_n(x,s,q)=xF_{n-1}(x,qs,q)+qsF_{n-2}(x,q^2s,q) \tag{3.2}$$

with initial values $F_0(x,s,q)=0$ and $F_1(x,s,q)=1.$

For consider the Fibonacci words which begin with a dot. The sum of their weights is $xF_{n-1}(x,qs,q).$ The weight of the set of all words which begin with a dash is $qsF_{n-2}(x,q^2s,q).$

If we consider the last elements of the Fibonacci words we get

$$F_n(x,s,q)=xF_{n-1}(x,s,q)+q^{n-2}sF_{n-2}(x,s,q). \tag{3.3}$$

This implies that

$$F_n(x,s,q)=\sum_{k=0}^{\left\lfloor\frac{n-1}{2}\right\rfloor} q^{k^2}\begin{bmatrix}n-1-k\\k\end{bmatrix}s^k x^{n-1-2k}, \tag{3.4}$$

where $\begin{bmatrix}n\\k\end{bmatrix}=\dfrac{[n]!}{[k]![n-k]!}$ denote the $q-$**binomial coefficients** which satisfy the recursions

$$\begin{bmatrix}n+1\\k\end{bmatrix}=q^k\begin{bmatrix}n\\k\end{bmatrix}+\begin{bmatrix}n\\k-1\end{bmatrix}=\begin{bmatrix}n\\k\end{bmatrix}+q^{n-k+1}\begin{bmatrix}n\\k-1\end{bmatrix}. \tag{3.5}$$

We prove (3.4) by induction. It is obviously true for $n=0$ and $n=1.$ The general case follows from

$$F_{n+1}(x,s,q)=xF_n(x,qs,q)+qsF_{n-1}(x,q^2s,q)=\sum_k q^{k^2}\begin{bmatrix}n-1-k\\k\end{bmatrix}(qs)^k x^{n-2k}+qs\sum_k q^{k^2}\begin{bmatrix}n-2-k\\k\end{bmatrix}\left(q^2s\right)^k x^{n-2-2k}$$

$$=\sum_k q^{k^2}\left(q^k\begin{bmatrix}n-1-k\\k\end{bmatrix}+\begin{bmatrix}n-1-k\\k-1\end{bmatrix}\right)s^k x^{n-2k}=\sum_k q^{k^2}\begin{bmatrix}n-k\\k\end{bmatrix}s^k x^{n-2k}.$$

As a special case we get that the weight of all Fibonacci words on $\{0,1,\cdots,n-2\}$ with $k$ dashes is $f(n,k,q)s^k x^{n-1-2k}$ with

$$f(n,k,q)=q^{k^2}\begin{bmatrix} n-1-k \\ k \end{bmatrix}. \tag{3.6}$$

The $q-$Fibonacci polynomials can be extended to negative values satisfying the same recurrence by

$$F_{-n}(x,s,q)=(-1)^{n-1}q^{\binom{n+1}{2}}\frac{F_n\left(x,\frac{s}{q^n},q\right)}{s^n}. \tag{3.7}$$

We note that $F_{-1}(x,s,q)=\frac{q}{s}$ and $F_{-2}(x,s,q)=-\frac{q^3x}{s^2}$.

**3.2.**

There is a useful analogue of the matrix powers.

Let $C(x,s)=\begin{pmatrix} 0 & 1 \\ s & x \end{pmatrix}$.

Then

$$M_n(x,s)=C(x,q^{n-1}s)C(x,q^{n-2}s)\cdots C(x,s)=\begin{pmatrix} sF_{n-1}(x,qs,q) & F_n(x,s,q) \\ sF_n(x,qs,q) & F_{n+1}(x,s,q)\end{pmatrix}. \tag{3.8}$$

Since $M_{m+n}(x,s)=M_m(x,q^ns)M_n(x,s)$ or

$$\begin{pmatrix} sF_{m+n-1}(x,qs,q) & F_{m+n}(x,s,q) \\ sF_{m+n}(x,qs,q) & F_{m+n+1}(x,s,q)\end{pmatrix}=\begin{pmatrix} q^nsF_{m-1}(x,q^{n+1}s,q) & F_m(x,q^ns,q) \\ q^nsF_m(x,q^{n+1}s,q) & F_{m+1}(x,q^ns,q)\end{pmatrix}\begin{pmatrix} sF_{n-1}(x,qs,q) & F_n(x,s,q) \\ sF_n(x,qs,q) & F_{n+1}(x,s,q)\end{pmatrix}$$

we get

$$F_{m+n}(x,s,q)=q^nsF_{m-1}(x,q^{n+1}s,q)F_n(x,s,q)+F_m(x,q^ns,q)F_{n+1}(x,s,q). \tag{3.9}$$

Such identities can of course also be derived from the combinatorial interpretation:

Consider all Fibonacci words on $\{0,1,\cdots,m+n-2\}$. If the number $n$ is covered by a dot or the beginning of a dash the weight of all words on $\{0,1,\cdots,n-1\}$ is $F_{n+1}(x,s,q)$ and the weight of all words on $\{n,\cdots m+n-2\}$ is $F_m(x,q^n s,q)$. Therefore the weight of all such words is $F_m(x,q^n s,q)F_{n+1}(x,s,q)$. If $\{n-1,n\}$ is covered by a dash then the weight of all words on $\{0,1,\cdots,n-2\}$ is $F_n(x,s,q)$ and the weight of the rest is $q^n sF_{m-1}(x,q^{n+1}s,q)$.

Since the trace of $M_n(x,s)$ is $F_{n+1}(x,s,q)+sF_{n-1}(x,qs,q)$ we get

**Theorem 3.2**

*Define the* $q-$***Lucas polynomials*** $L_n(x,s,q)$ *by*

$$L_n(x,s,q)=tr\left(M_n(x,s)\right). \tag{3.10}$$

*Then*

$$L_n(x,s,q)=F_{n+1}(x,s,q)+sF_{n-1}(x,qs,q). \tag{3.11}$$

They have the same initial values as for $q=1$.

For $n=0$ the matrix $C(x,q^{n-1}s)C(x,q^{n-2}s)\cdots C(x,s)$ reduces to the identity matrix which gives $L_0(x,s,q)=2$. The right-hand side reduces also to $2$ since $F_1(x,s,q)+sF_{-1}(x,qs,q)=1+s\frac{q}{qs}=2$.

For $n=1$ we get $L_1(x,s,q)=F_2(x,s,q)+sF_0(x,qs,q)=x$.

The determinant of (3.8) gives a $q-$Cassini formula

$$F_{n-1}(x,qs,q)F_{n+1}(x,s,q)-F_n(x,s,q)F_n(x,qs,q)=(-1)^n q^{\binom{n}{2}}s^{n-1}. \tag{3.12}$$

**3.3.**

We can also interpret the $q-$Lucas polynomials as the weight of all circular Fibonacci words on $\{0,1,\cdots,n-1\}$.

If the arc $0$ is covered by a dot or the first point of a dash then the weight of all such Fibonacci words is $F_{n+1}(x,s,q)$. If $0$ is the endpoint of a dash then this dash has weight $s$ and the remaining word has weight $F_{n-1}(x,qs,q)$.

Therefore the weight of all coverings is $L_n(x,s,q)=F_{n+1}(x,s,q)+sF_{n-1}(x,qs,q)$.

Let $L_n(x,s,q)=\sum_{k=0}^{\left\lfloor \frac{n}{2}\right\rfloor} l(n,k,q)s^k x^{n-2k}$.

Then we get from (3.11)

$$l(n,k,q)= f(n+1,k,q)+q^{k-1}f(n-1,k-1,q)=q^{k^2-k}\frac{[n]}{[n-k]}\begin{bmatrix} n-k \\ k \end{bmatrix}. \tag{3.13}$$

Thus we have for $n>0$

$$L_n(x,s,q)=\sum_{k=0}^{\left\lfloor \frac{n}{2}\right\rfloor} q^{k^2-k}\frac{[n]}{[n-k]}\begin{bmatrix} n-k \\ k \end{bmatrix} s^k x^{n-2k} \tag{3.14}$$

with initial value $L_0(x,s,q)=2$.

Another formula for the $q-$Lucas polynomials is

$$xL_n(x,qs,q)=F_{n+2}\left(x,s,q\right)-q^{n+1}s^2F_{n-2}(x,q^2s,q). \tag{3.15}$$

For

$$xL_n(x,qs,q)=xF_{n+1}(x,qs,q)+qsxF_{n-1}(x,q^2s,q)=F_{n+2}(x,s,q)-qsF_n(x,q^2s,q)+qsxF_{n-1}(x,q^2s,q)$$
$$=F_{n+2}(x,s,q)-qs\left(F_n(x,q^2s,q)-xF_{n-1}(x,q^2s,q)\right)=F_{n+2}(x,s,q)-q^{n+1}s^2F_{n-2}(x,q^2s,q).$$

Comparing coefficients in (3.15) we see that

$$q^k l(n,k,q)=q^{k^2}\begin{bmatrix} n-k+1 \\ k \end{bmatrix}-q^{n+(k-1)^2}\begin{bmatrix} n-k-1 \\ k-2 \end{bmatrix}. \tag{3.16}$$

It is easily verified that (3.11) gives

$$L_{-n}(x,s,q)=(-1)^n\frac{q^{\binom{n+1}{2}}}{s^n}L_n\left(x,\frac{s}{q^n},q\right). \tag{3.17}$$

Contrary to the $q-$Fibonacci polynomials the $q-$Lucas polynomials do not satisfy a simple recurrence. The reason is that $F_{n+1}(x,s,q)$ and $sF_{n-1}(x,qs,q)$ do not satisfy the same recurrence. For $H_n(s)=F_{n+1}(x,s,q)$ satisfies $H_n(s)=xH_{n-1}(qs)+qsH_{n-2}(q^2s)$, whereas $K_n(s)=sF_{n-1}(x,qs,q)$ satisfies $K_n(s)=q^{-1}xK_{n-1}(qs)+sK_{n-2}(q^2s).$

The $q-$Zeilberger algorithm gives a computer proof of

**Theorem 3.3**

*The $q-$Lucas polynomials satisfy the recurrence relation*

$$\begin{aligned}L_n(x,s,q)&=(1+q)xL_{n-1}(x,s,q)+\left(q^{n-2}(1+q)s-qx^2\right)L_{n-2}(x,s,q)\\&-q^{n-2}(1+q)xsL_{n-3}(x,s,q)-q^{2n-5}s^2L_{n-4}(x,s,q).\end{aligned} \tag{3.18}$$

A "human" proof runs as follows:

Let $H_n(s)$ satisfy

$$H_n(s)-(1+q)xH_{n-1}(s)-\left(q^{n-2}(1+q)s-qx^2\right)H_{n-2}(s)+q^{n-2}(1+q)xsH_{n-3}(s)+q^{2n-5}s^2H_{n-4}(s)=0.$$

Then it is easily verified that $sH_{n-1}(qs)$ satisfies the same recurrence.

Because of (3.11) it suffices to show that

$$\begin{aligned}&F_n(x,s,q)-(1+q)xF_{n-1}(x,s,q)-\left(q^{n-2}(1+q)s-qx^2\right)F_{n-2}(x,s,q)\\&+q^{n-2}(1+q)xsF_{n-3}(x,s,q)+q^{2n-5}s^2F_{n-4}(x,s,q)=0.\end{aligned}$$

This can be reduced to several applications of the recursion (3.3). We get successively

$$\begin{aligned}0&=qxF_{n-1}(x,s,q)+q^{n-1}sF_{n-2}(x,s,q)-qx^2F_{n-2}(x,s,q)\\&-q^{n-2}xsF_{n-3}(x,s,q)-q^{n-1}xsF_{n-3}(x,s,q)-q^{2n-5}s^2F_{n-4}(x,s,q),\end{aligned}$$

$$0=F_n(x,s,q)-xF_{n-1}(x,s,q)-q^{n-2}xsF_{n-3}(x,s,q)-q^{2n-6}s^2F_{n-4}(x,s,q)$$

and finally $0=F_n(x,s,q)-xF_{n-1}(x,s,q)-q^{n-2}sF_{n-2}(x,s,q).$

## Remarks

As far as I know the first one who studied $q-$analogues of the Fibonacci numbers was I. Schur [22]. In 1917 he introduced $F_n(1,1,q)$ and $F_n(1,q,q)$ as certain determinants and proved his celebrated polynomial versions of the **Rogers-Ramanujan identities**

$$F_{n+1}(1,1,q)=\sum_{k\in\mathbb{Z}}(-1)^k q^{\frac{k(5k-1)}{2}}\begin{bmatrix} n \\ \left\lfloor \frac{n+5k}{2} \right\rfloor \end{bmatrix} \tag{3.19}$$

and

$$F_n(1,q,q)=\sum_{k\in\mathbb{Z}}(-1)^k q^{\frac{k(5k-3)}{2}}\begin{bmatrix} n \\ \left\lfloor \frac{n+5k-1}{2} \right\rfloor \end{bmatrix}. \tag{3.20}$$

In [4] and [5] L. Carlitz systematically studied $q-$analogues of the Fibonacci and Lucas numbers and of Fibonacci polynomials. Therefore I have associated these polynomials with his name although he seemingly never defined $q-$Lucas polynomials.

His starting point was the observation that the number of sequences of zeros and ones $(a(1),a(2),\cdots,a(n))$ of length $n$ in which consecutive $1$'s are forbidden is equal to the Fibonacci number $F_{n+2}$. In our model this amounts to choosing $a(i)=1$ if $i$ is the endpoint of a dash and $a(i)=0$ in all other cases. The sequences obtained in this way are characterized by $a(0)=0$ and the fact that no consecutive ones occur. Carlitz used this observation to define a $q-$analogue of the Fibonacci number $F_n(q)$ by $\sum q^{a(1)+2a(2)+\cdots+(n-2)a(n-2)}$ where the sum is extended over all $(a(1),a(2),\cdots,a(n-2))$ of zeros and ones where consecutive ones are forbidden. In our notation $F_n(q)=F_n(1,1,q)$.

His $q-$analogue of the Lucas numbers ([4], (1.8)) is in our notation $L_n(1,q,q)=F_{n+2}(1,1,q)-q^{n+1}F_{n-2}(1,q^2,q)$.

Note that

$$l(n,k,q)=\sum_{a(i)} q^{\sum_{i=0}^{n-1} ia(i)}, \tag{3.21}$$

where the sum extends over all $(a(0),a(1),\cdots,a(n-1))$ of $a(i)\in\{0,1\}$ such that consecutive ones and also $a(0)=a(n-1)=1$ are not allowed.

Carlitz gave a direct proof that (3.16) coincides with $\Sigma' q^{a(1)+2a(2)+\cdots+na(n)}$ where $\Sigma'$ runs over all sequences $(a(1),\cdots,a(n))$ with no consecutive ones and where $a(1)=a(n)=1$ is forbidden.

It seems that there are some properties of the classical Lucas polynomials which cannot be generalized to $q-$Lucas polynomials. For example there is no analogue of (2.21) since (3.9) for $m=n$ has no factorization.

Another instance where Lucas polynomials occur is in the recurrence relation for subsequences $F_{\ell n}(x,s)$. They satisfy $F_{\ell n}(x,s)-L_{\ell}(x,s)F_{\ell(n-1)}(x,s)+(-s)^{\ell}F_{\ell(n-2)}(x,s)$. But the corresponding $q-$Fibonacci sequences satisfy instead (cf. [13])

$$F_{\ell n}(x,s,q)-\frac{F_{2\ell}(x,s,q)}{F_{\ell}(x,q^{\ell}s,q)}F_{\ell(n-1)}(x,q^{\ell}s,q)+(-s)^{\ell}q^{\frac{\ell(3\ell-1)}{2}}\frac{F_{\ell}(x,s,q)}{F_{\ell}(x,q^{\ell}s,q)}F_{\ell(n-2)}(x,q^{2\ell}s,q).$$

The coefficients are in general not even polynomials.

Several papers (cf. [1], [2], [3], [7], [8], [9], [10],[11], [12], [13], [16], [17], [18],[19],[21], [22]) deal with combinatorial interpretations of these polynomials, give applications to the Rogers-Ramanujan formulas or prove $q-$analogues of some of the almost inexhaustible set of identities satisfied by Fibonacci and Lucas numbers or polynomials.

## 4. Another class of q-Fibonacci and q- Lucas polynomials

### 4.1.

To obtain another class of $q-$Fibonacci and $q-$Lucas polynomials we consider words $c=c_1c_2\cdots c_m$ of letters $c_i\in\{a,b\}$ and associate with $c$ the **weight**

$$W(c)=W(c)(s)=q^{i_1+i_2+\cdots+i_k}s^kx^{m-k}, \tag{4.1}$$

if $c_{i_1}=c_{i_2}=\cdots=c_{i_k}=b, 1\le i_1<\cdots<i_k\le m,$ and all other $c_i=a.$ The weight of the empty word $\varepsilon$ is defined to be $W(\varepsilon)=1.$

We then have

$$\begin{aligned} W(ac)(s)&=xW(c)(qs),\\ W(bc)(s)&=qsW(c)(qs),\\ W(ca)(s)&=xW(c)(s),\\ W(cb)(s)&=q^{m+1}sW(c)(s). \end{aligned} \tag{4.2}$$

**Definition 4.1**

*The* $q-$ ***Fibonacci polynomial*** $Fib_n(x,s,q)$ *is defined by*

$$Fib_n(x,s,q)=\sum_c W(c) \tag{4.3}$$

*where the sum is extended over all Fibonacci words of length* $n-1$.

By considering the first letter of each word we see from (4.2) that

$$Fib_n(x,s,q)=xFib_{n-1}(x,qs,q)+qsFib_{n-2}(x,qs,q). \tag{4.4}$$

The initial values are given by $Fib_0(x,s,q)=0$ and $Fib_1(x,s,q)=1$.

By considering the last letter of each word we get in the same way

$$Fib_n(x,s,q)=xFib_{n-1}(x,s,q)+q^{n-2}sFib_{n-2}\left(x,\frac{s}{q},q\right). \tag{4.5}$$

To obtain the second term let us suppose that the Fibonacci sequence $cb$ of order $n$ has $k$ letters $b$. Then $W(cb)(s)=\left(q^{n-k-1}s\right)q^{i_1+\cdots+i_{k-1}}s^{k-1}x^{n-k-1}=\left(q^{n-2}s\right)q^{i_1+\cdots+i_{k-1}}\left(\frac{s}{q}\right)^{k-1}x^{n-k-1}=q^{n-2}sW(c)\left(\frac{s}{q}\right)$.

Since this expression is independent of $k$, we get the second term.

The first polynomials are

$0,1,x,x^2+qs,x^3+(1+q)qsx,x^4+qs[3]x^2+q^3s^2,\cdots$

**Theorem 4.2**

$$Fib_n(x,s,q)=\sum_{k=0}^{\left\lfloor\frac{n-1}{2}\right\rfloor} q^{\binom{k+1}{2}}\begin{bmatrix} n-1-k \\ k\end{bmatrix}s^k x^{n-1-2k}. \tag{4.6}$$

This can be proved by induction using (4.4).

**4.2.**

These polynomials have a very special property as has been observed in [6] and [8]. They arise by replacing in the classical Fibonacci polynomials the variable $x$ by the operator $x+(q-1)sD$ and applying these operators to the constant polynomial $1$.

**Theorem 4.3**

*Let* $D$ *be the* $q-$*differentiation operator defined by* $Df(x)=\frac{f(x)-f(qx)}{(1-q)x}$. *Then* $Fib_n(x,s,q)$ *satisfies*

$$Fib_n(x,s,q)=F_n\left(x+(q-1)sD,s\right)1. \qquad (4.7)$$

*With other words this means that*

$$Fib_n(x,s,q)=xFib_{n-1}(x,s,q)+(q-1)sDFib_{n-1}(x,s,q)+sFib_{n-2}(x,s,q). \qquad (4.8)$$

We call a recurrence where the operator $(q-1)sD$ occurs a $D-$**recurrence**.

Comparing coefficients of $s^k$ in $Fib_{n+1}(x,s,q)$ this amounts to

$$q^{\binom{k+1}{2}}\begin{bmatrix}n-k\\k\end{bmatrix}=q^{\binom{k+1}{2}}\begin{bmatrix}n-1-k\\k\end{bmatrix}+\left(q^{n+1-2k}-1\right)q^{\binom{k}{2}}\begin{bmatrix}n-k\\k-1\end{bmatrix}+q^{\binom{k}{2}}\begin{bmatrix}n-1-k\\k-1\end{bmatrix}.$$

This is equivalent with

$$q^k\left(\begin{bmatrix}n-k\\k\end{bmatrix}-\begin{bmatrix}n-1-k\\k\end{bmatrix}\right)=\left(q^{n+1-2k}-1\right)\begin{bmatrix}n-k\\k-1\end{bmatrix}+\begin{bmatrix}n-1-k\\k-1\end{bmatrix}$$

or

$$\left(q^{n-k}-1\right)\begin{bmatrix}n-1-k\\k-1\end{bmatrix}=\left(q^{n+1-2k}-1\right)\begin{bmatrix}n-k\\k-1\end{bmatrix}$$

which is obviously true.

As an example consider $F_4(x,s)=x^3+2sx.$

We compute $F_4(x+(q-1)sD,s)1=(x+(q-1)sD)^3 1+2s(x+(q-1)sD)1.$

$$\left(x+(q-1)sD)\right)1 = x,$$
$$\left(x+(q-1)sD)\right)x = x^2+(q-1)s,$$
$$\left(x+(q-1)sD)\right)\left(x^2+(q-1)s\right) = x^3+(q+q^2-2)sx$$

Therefore $F_4(x+(q-1)sD,s)1 = x^3+(q+q^2-2)sx+2sx = x^3+(q+q^2)sx.$

**4.3.**

We define the corresponding $q-$**Lucas polynomials** by

$$Luc_n(x,s,q) = L_n\left(x+(q-1)sD,s\right)1. \tag{4.9}$$

The first polynomials are

$2, x, x^2+(1+q)s, x^3+[3]sx, x^4+[4]sx^2+q(1+q^2)s^2, \cdots$

By applying the linear map

$$f(x) \to f(x+(q-1)sD)1 \tag{4.10}$$

to (2.18) we get

$$Luc_n(x,s,q) = Fib_{n+1}(x,s,q)+sFib_{n-1}(x,s,q) \tag{4.11}$$

for $n>0$.

This implies the explicit formula

$$Luc_n(x,s,q) = \sum_{k=0}^{\frac{n}{2}} q^{\binom{k}{2}} \frac{[n]}{[n-k]} \begin{bmatrix} n-k \\ k \end{bmatrix} s^k x^{n-2k} \tag{4.12}$$

for $n>0$.

For the proof observe that

$$q^k \begin{bmatrix} n-k \\ k \end{bmatrix} + \begin{bmatrix} n-k-1 \\ k-1 \end{bmatrix} = \frac{q^k[n-k]+[k]}{[n-k]} \begin{bmatrix} n-k \\ k \end{bmatrix} = \frac{[n]}{[n-k]} \begin{bmatrix} n-k \\ k \end{bmatrix}.$$

(4.11) has the following combinatorial interpretation:

Consider all circular Fibonacci words on $\{0,1,\cdots,n-1\}.$

If the arc $0$ is covered by a dot or the first point of a dash then the weight of all such words is $Fib_{n+1}(x,s,q).$ If $0$ is the endpoint of a dash we split this $b$ into $b=b_0b_1$ and associate with this covering the word $b_1c_1\cdots c_mb_0$ and define its weight as $sW(c_1\cdots c_m).$

**4.4.**

For $Fib_n(x,s,q)$ and $Luc_n(x,s,q)$ there is a weak **analogue of Binet's formulas**. Let $A$ be the linear operator $A=x+(q-1)sD$ on the polynomials.

Define the formal expressions

$$\alpha(q)=\frac{A+\sqrt{A^2+4s}}{2} \tag{4.13}$$

and

$$\beta(q)=\frac{A-\sqrt{A^2+4s}}{2}. \tag{4.14}$$

We call these "formal expressions" because we do not give an interpretation for $\sqrt{A^2+4s}$. In the formulas below these square roots cancel.

Since $\alpha(q)^2-A\alpha(q)-s=\beta(q)^2-A\beta(q)-s=0$ the sequences $\left(\alpha(q)^n\right)_{-\infty}^{\infty}$ and $\left(\beta(q)^n\right)_{-\infty}^{\infty}$ satisfy the recurrence

$$\alpha(q)^n-A\alpha(q)^{n-1}-s\alpha(q)^{n-2}=\beta(q)^n-A\beta(q)^{n-1}-s\beta(q)^{n-2}=0$$

for all $n\in\mathbb{Z}.$

Since the $q-$Fibonacci and the $q-$Lucas polynomials satisfy the same recurrence we get from the initial values

$$Luc_n(x,s,q)=\left(\alpha(q)^n+\beta(q)^n\right)1 \tag{4.15}$$

and

$$Fib_n(x,s,q)=\frac{\alpha(q)^n-\beta(q)^n}{\alpha(q)-\beta(q)}1 \tag{4.16}$$

for $n \geq 0$. Note that in these expressions the term $\sqrt{A^2+4s}$ does not occur. Therefore these are well-defined polynomials.

We can use these identities to extend these polynomials to negative $n$.

We then get for $n > 0$

$$Luc_{-n}(x,s,q) = \left(\alpha(q)^{-n} + \beta(q)^{-n}\right)1 = (-1)^n \frac{\beta(q)^n + \alpha(q)^n}{s^n} = (-1)^n \frac{Luc_n(x,s,q)}{s^n} \tag{4.17}$$

and

$$Fib_{-n}(x,s,q) = \frac{\alpha(q)^{-n} - \beta(q)^{-n}}{\alpha(q) - \beta(q)}1 = (-1)^{n-1} \frac{Fib_n(x,s,q)}{s^n}. \tag{4.18}$$

**4.5.**

If we consider $C(A,s) = \begin{pmatrix} 0 & 1 \\ s & A \end{pmatrix}$ then we get

$$C(A,s)^n 1 = \begin{pmatrix} sFib_{n-1}(x,s,q) & Fib_n(x,s,q) \\ sFib_n(x,s,q) & Fib_{n+1}(x,s,q) \end{pmatrix}. \tag{4.19}$$

So we also have

$$Luc_n(x,s,q) = tr\left(C(A,s)^n 1\right). \tag{4.20}$$

These $q-$Lucas polynomials satisfy also the recurrence

$$Luc_n(x,s,q) = xLuc_{n-1}(x,s,q) + (q-1)sDLuc_{n-1}(x,s,q) + sLuc_{n-2}(x,s,q). \tag{4.21}$$

But the $q-$Lucas polynomials do not satisfy the recurrence (4.4).

We have instead

$$Luc_n(x,s,q) - xLuc_{n-1}(x,qs,q) - qsLuc_{n-2}(x,qs,q) = (1-q)sFib_{n-1}(x,s,q). \tag{4.22}$$

For replacing $n \to n+1$ and since $Fib_n(x,s,q)$ satisfies the recurrence this is equivalent with

$$sFib_n(x,s,q) - xqsFib_{n-1}(x,qs,q) - qsqsFib_{n-2}(x,qs,q) = (1-q)sFib_n(x,s,q)$$

or

$$Fib_n(x,s,q) - xFib_{n-1}(x,qs,q) - qsFib_{n-2}(x,qs,q) = 0.$$

**4.6.**

Till now we have no recurrences for these new $q-$polynomials with fixed values $x$ and $s$.

For the polynomials $Fib_n(x,s,q)$ such a recurrence can be found by combining (4.4) and (4.5).

This gives

$$Fib_n(x,s,q) = xFib_{n-1}(x,s,q) + q^{n-2}sxFib_{n-3}(x,s,q) + q^{n-2}s^2Fib_{n-4}(x,s,q). \qquad (4.23)$$

For the $q-$Lucas polynomials we have a rather ugly recursion (cf. [14]).

$$Luc_{n+4}(x,s,q) = xLuc_{n+3}(x,s,q) - q^{n+1}\frac{[2]}{[n+1]}Luc_{n+2}(x,s,q)$$
$$+q^{n+1}\frac{[n+3]}{[n+1]}sxLuc_{n+1}(x,s,q) + q^{n+1}\frac{[n+3]}{[n+1]}s^2Luc_n(x,s,q).$$

**4.7.**

We define now the modified $q-$**Lucas polynomials** by $Luc_n^*(x,s,q) = Luc_n(x,s,q)$ for $n > 0$ and $Luc_0^*(x,s,q) = 1.$

A $q-$analogue of (2.24) with applications to Rogers-Ramanujan type identities (cf. [14]) is

**Theorem 4.4**

$$\sum_{k=0}^{\left\lfloor \frac{n}{2} \right\rfloor} (-s)^k \begin{bmatrix} n \\ k \end{bmatrix} Luc_{n-2k}^*(x,s,q) = x^n. \qquad (4.24)$$

**Proof**

This is trivially true for $n = 0$ and $n = 1.$ We set $Luc_n^*(x,s,q) = 0$ for $n < 0$ and get by induction

$$\sum_k (-s)^k \begin{bmatrix} n+1 \\ k \end{bmatrix} Luc^*_{n+1-2k}(x,s,q) = \sum_k (-s)^k \left( q^k \begin{bmatrix} n \\ k \end{bmatrix} + \begin{bmatrix} n \\ k-1 \end{bmatrix} \right) Luc^*_{n+1-2k}(x,s,q)$$
$$= \sum_k (-s)^k \left( \begin{bmatrix} n \\ k \end{bmatrix} + \begin{bmatrix} n \\ k-1 \end{bmatrix} \right) Luc^*_{n+1-2k}(x,s,q) + \sum_k (-s)^k \left( \left(q^k - 1\right) \begin{bmatrix} n \\ k \end{bmatrix} \right) Luc^*_{n+1-2k}(x,s,q)$$
$$= \sum_k (-s)^k \begin{bmatrix} n \\ k \end{bmatrix} \left( Luc^*_{n+1-2k}(x,s,q) - sL^*_{n-1-2k}(x,s,q) \right) + \sum_k (-s)^k \left( \left(q^n - 1\right) \begin{bmatrix} n-1 \\ k-1 \end{bmatrix} \right) Luc^*_{n+1-2k}(x,s,q)$$
$$= \sum_{k=0} (-s)^k \begin{bmatrix} n \\ k \end{bmatrix} \left( x + (q-1)sD \right) \left( Luc^*_{n-2k}(x,s,q) \right) - \left(q^n - 1\right) \sum_k (-s)^k \left( \begin{bmatrix} n-1 \\ k \end{bmatrix} \right) Luc^*_{n-1-2k}(x,s,q)$$
$$= \left( x + (q-1)sD \right) x^n - \left(q^n - 1\right) x^{n-1} = x^{n+1}.$$

We have only to observe that

$$Luc^*_{n+1}(x,s,q) - sLuc^*_{n-1}(x,s,q) = \left( x + (q-1)sD \right) \left( Luc^*_n(x,s,q) \right)$$

for $n \neq 1$ and that

$$(-s)^m \begin{bmatrix} 2m+1 \\ m \end{bmatrix} \left( Luc^*_2(x,s) - sLuc^*_0(x,s) \right) + (-s)^{m+1} \begin{bmatrix} 2m+1 \\ m+1 \end{bmatrix} = (-s)^m \begin{bmatrix} 2m+1 \\ m \end{bmatrix} (x^2 + qs - s)$$
$$= \left( x + (q-1)sD \right) x.$$

An interesting $q-$analogue of $F_n(1,-1)$ is $Fib_n\left(1, -\frac{1}{q}, q\right)$.

Here it is easily verified (cf. [8]) that

$$Fib_{3n}\left(1, -\frac{1}{q}, q)\right) = 0, Fib_{3n+1}\left(1, -\frac{1}{q}, q)\right) = (-1)^n q^{\frac{n(3n-1)}{2}}, Fib_{3n+2}\left(1, -\frac{1}{q}, q)\right) = (-1)^n q^{\frac{n(3n+1)}{2}}. \qquad (4.25)$$

As shown in [14] these polynomials give an easy approach to **Slater's Bailey pairs** A(1)-A(8).

## Remarks

These polynomials have been introduced in [6] and [8]. They also occur in [14],[15],[20].

As far as I know $D-$recurrences have not been studied in other contexts.

In [15] we have defined a new $q-$analogue of the **Hermite polynomials** $H_n(x,s \mid q) = (x - sD)^n 1.$

By applying the linear map (4.10) to the identity (2.24) we get

$$H_n(x,(q-1)s\,|\,q)=\sum_{k=0}^{\left\lfloor\frac{n}{2}\right\rfloor}\binom{n}{k}s^k Luc^*_{n-2k}(x,-s,q). \tag{4.26}$$

## 5. A common generalization

### 5.1.

By comparing coefficients in (3.4) and (4.6) we see that the following result holds.

**Theorem 5.1**

*Let $U$ be the linear isomorphism on the polynomials in $s$ defined by*

$$Us^k=q^{\binom{k}{2}}s^k. \tag{5.1}$$

*Then for each $j\in\mathbb{Z}$*

$$U\left(s^j Fib_n(x,q^i s,q)\right)=s^j q^{\binom{j}{2}}F_n(x,q^{i+j}s,q). \tag{5.2}$$

This can be used to translate identities for one class of $q-$Fibonacci polynomials to the other class. Consider for example the identity (cf. [8])

$$Fib_{2n}(x,s,q)=\sum_{k=0}^{n}\begin{bmatrix}n\\k\end{bmatrix}q^{\binom{k+1}{2}}s^k x^{n-k}Fib_{n-k}(x,q^n s,q). \tag{5.3}$$

If we apply $U$ we get

$$F_{2n}(x,s,q)=\sum_{k=0}^{n}\begin{bmatrix}n\\k\end{bmatrix}q^{k^2}s^k x^{n-k}F_{n-k}(x,q^{n+k}s,q). \tag{5.4}$$

In the same way the identity (cf. [14])

$$\sum_{k=0}^{n}(-1)^k q^{\binom{k}{2}}\begin{bmatrix}n\\k\end{bmatrix}x^k Fib_{2n+m-k}(x,s,q)=q^{\binom{n}{2}+mn}s^n Fib_m\left(x,\frac{s}{q^n},q\right) \tag{5.5}$$

is transformed into

$$\sum_{k=0}^{n}(-1)^k q^{\binom{k}{2}}\begin{bmatrix}n\\k\end{bmatrix}x^k F_{2n+m-k}(x,s,q)=q^{2\binom{n}{2}+mn}s^n F_m(x,s,q). \tag{5.6}$$

It is also instructive to see how the polynomials with negative indices are transformed.

$$U\left(Fib_{-n}(x,s,q)\right)=U\left((-1)^{n-1}s^{-n}Fib_n(x,s,q)\right)=(-1)^{n-1}s^{-n}q^{\binom{-n}{2}}F_n\left(x,\frac{s}{q^n},q\right)=F_{-n}(x,s,q).$$

Theorem 5.1 implies that also for the Carlitz $q-$Fibonacci and $q-$Lucas polynomials an analogous $D-$recurrence holds.

**Corollary 5.2**

*Let* $D$ *be the* $q-$*differentiation operator on the polynomials in* $x$. *Then*

$$F_n(x,s,q)=xF_{n-1}(x,s,q)+(q-1)sDF_{n-1}(x,qs,q)+sF_{n-2}(x,qs,q) \tag{5.7}$$

*and*

$$L_n(x,s,q)=xL_{n-1}(x,s,q)+(q-1)sDL_{n-1}(x,qs,q)+sL_{n-2}(x,qs,q). \tag{5.8}$$

*More generally this* $D-$*recurrence holds for all linear combinations of polynomials* $s^kF_n(x,q^ks,q)$ *with complex coefficients.*

**Proof**

$$F_n(x,s,q)=UFib_n(x,s,q)=xUFib_{n-1}(x,s,q)+(q-1)DUsFib_{n-1}(x,s,q)+UsFib_{n-2}(x,s,q)$$
$$=xF_{n-1}(x,s,q)+(q-1)sDF_{n-1}(x,qs,q)+sF_{n-2}(x,qs,q).$$

If we apply $U$ to $s^kFib_n(x,s,q)$ we get

$$s^kq^{\binom{k}{2}}F_n(x,q^ks,q)=xs^kq^{\binom{k}{2}}F_{n-1}(x,q^ks,q)+q^{\binom{k+1}{2}}(q-1)Ds^{k+1}F_{n-1}(x,q^{k+1}s,q)+q^{\binom{k+1}{2}}s^{k+1}F_{n-2}(x,q^{k+1}s,q).$$

Dividing by $q^{\binom{k}{2}}$ and setting $G_n(s)=s^kF_n(x,q^ks,q)$ this gives

$$G_n(s)=xG_{n-1}(s)+(q-1)sDG_{n-1}(qs)+sG_{n-2}(qs).$$

**5.2.**

These results suggest the following

**Theorem 5.3**

*Let $U_m$ be the linear operator on the polynomials in $s$ defined by*

$$U_m s^k = q^{m\binom{k}{2}} s^k \tag{5.9}$$

*and define polynomials $\Phi_n(x,s,m,q)$ by*

$$\Phi_n(x,s,m,q) := U_m Fib_n(x,s,q), \tag{5.10}$$

i.e.

$$\Phi_n(x,s,m,q) = \sum_{k=0}^{\left\lfloor \frac{n-1}{2} \right\rfloor} q^{\binom{k+1}{2}+m\binom{k}{2}} \begin{bmatrix} n-1-k \\ k \end{bmatrix} s^k x^{n-1-2k}. \tag{5.11}$$

*Then for each $j \in \mathbb{Z}$*

$$U_m\left(s^j Fib_n(x,q^i s,q)\right) = s^j q^{m\binom{j}{2}} \Phi_n(x,q^{i+mj}s,m,q). \tag{5.12}$$

*These polynomials satisfy the $D-$recurrrence*

$$\Phi_n(x,s,m,q) = x\Phi_{n-1}(x,s,m,q) + (q-1)sD\Phi_{n-1}(x,q^m s,m,q) + s\Phi_{n-2}(x,q^m s,m,q) \tag{5.13}$$

*with initial values $\Phi_0(x,s,m,q) = 0$ and $\Phi_1(x,s,m,q) = 1$.*

*More generally this $D-$recurrence holds for all linear combinations of polynomials $s^k\Phi_n(x,q^{km}s,q)$ with complex coefficients.*

**Proof**

We first prove (5.12)

$$U_m\left(s^j Fib_n(x,q^i s,q)\right)=U_m\left(\sum_{k=0}^{\left\lfloor\frac{n-1}{2}\right\rfloor} q^{\binom{k+1}{2}}\begin{bmatrix}n-1-k\\k\end{bmatrix} q^{ik}s^{k+j}x^{n-1-2k}\right)$$

$$=\sum_{k=0}^{\left\lfloor\frac{n-1}{2}\right\rfloor} q^{\binom{k+1}{2}}\begin{bmatrix}n-1-k\\k\end{bmatrix} q^{ik}q^{m\binom{k+j}{2}}s^{k+j}x^{n-1-2k}=s^j q^{m\binom{j}{2}}\Phi_n(x,q^{i+jm}s,m,q).$$

Applying $U_m$ to (4.8)

gives

$$\Phi_n(x,s,m,q)=U_m Fib_n(x,s,q)=x\Phi_{n-1}(x,s,m,q)+(q-1)Ds\Phi_{n-1}(x,q^m s,m,q)+s\Phi_{n-2}(x,q^m s,m,q).$$

Finally if we apply $U_m$ to $s^k Fib_n(x,s,q)$ we get

$$s^k q^{m\binom{k}{2}}\Phi_n(x,q^{mk}s,q)=xs^k q^{m\binom{k}{2}}\Phi_{n-1}(x,q^{mk}s,q)$$
$$+q^{\binom{k+1}{2}m}(q-1)Ds^{k+1}\Phi_{n-1}(x,q^{(k+1)m}s,q)+q^{\binom{k+1}{2}m}s^{k+1}\Phi_{n-2}(x,q^{(k+1)m}s,q).$$

Dividing by $q^{m\binom{k}{2}}$ and setting $G_n(s)=s^k\Phi_n(x,q^{km}s,q)$ this gives

$$G_n(s)=xG_{n-1}(s)+(q-1)sDG_{n-1}(q^m s)+sG_{n-2}(q^m s).$$

We define corresponding $q-$Lucas polynomials by

$$\Lambda_n(x,s,m,q)=U_m Luc_n(x,s,q)=\Phi_{n+1}(x,s,m,q)+s\Phi_{n-1}(x,q^m s,m,q). \tag{5.14}$$

They are given by

$$\Lambda_n(x,s,m,q)=\sum_{k=0}^{\frac{n}{2}} q^{(m+1)\binom{k}{2}}\frac{[n]}{[n-k]}\begin{bmatrix}n-k\\k\end{bmatrix}s^k x^{n-2k}. \tag{5.15}$$

By Theorem 5.3 they also satisfy

$$\Lambda_n(x,s,m,q) = x\Lambda_{n-1}(x,s,m,q) + (q-1)sD\Lambda_{n-1}(x,q^m s,m,q) + s\Lambda_{n-2}(x,q^m s,m,q). \quad (5.16)$$

**5.3.**

From (4.4) we conclude the recurrence

$$\Phi_n(x,s,m,q) = x\Phi_{n-1}(x,qs,m,q) + qs\Phi_{n-2}(x,q^{m+1}s,m,q) \quad (5.17)$$

and from (4.5)

$$\Phi_n(x,s,m,q) = x\Phi_{n-1}(x,s,m,q) + q^{n-2}s\Phi_{n-2}(x,q^{m-1}s,m,q). \quad (5.18)$$

In order to translate this into algebraic terms we define the linear operator $\eta$ on the polynomials in the indeterminate $s$ by

$$\eta f(s) = f(qs). \quad (5.19)$$

Then (5.17) can be written in the form

$$\Phi_n(x,s,m,q) = x\eta\Phi_{n-1}(x,s,m,q) + qs\eta^{m+1}\Phi_{n-2}(x,s,m,q). \quad (5.20)$$

Thus we are led to consider the ring homomorphism $\varphi_m : P(a,b) \to \mathbb{C}(q)[x\eta, qs\eta^{m+1}]$ from the ring of the non-commutative polynomials in $a,b$ to the (also non-commutative) ring of polynomials in the operators $x\eta$ and $qs\eta^2$. These operators satisfy $(x\eta)\left(qs\eta^{m+1}\right) = q^2xs\eta^{m+2} = q\left(qs\eta^{m+1}\right)(x\eta)$.

This implies that

$$\varphi_m\left((a+b)^n\right) = \left(qs\eta^{m+1} + x\eta\right)^n = \sum_{k=0}^{n} \begin{bmatrix} n \\ k \end{bmatrix} \left(qs\eta^{m+1}\right)^k (x\eta)^{n-k} = \sum_{k=0}^{n} q^{k^2} \begin{bmatrix} n \\ k \end{bmatrix} s^k x^{n-k} \eta^{n+mk}. \quad (5.21)$$

This is easily proved by induction:

$$\left(qs\eta^{m+1}+x\eta\right)^{n+1}=\left(qs\eta^{m+1}+x\eta\right)\sum_k \begin{bmatrix} n \\ k \end{bmatrix}\left(qs\eta^{m+1}\right)^k (x\eta)^{n-k}=\sum_k \begin{bmatrix} n \\ k \end{bmatrix}\left(qs\eta^{m+1}\right)^{k+1}(x\eta)^{n-k}$$

$$+\sum_k \begin{bmatrix} n \\ k \end{bmatrix} x\eta\left(qs\eta^{m+1}\right)^k (x\eta)^{n-k}$$

$$=\sum_k \begin{bmatrix} n \\ k-1 \end{bmatrix}\left(qs\eta^{m+1}\right)^k (x\eta)^{n+1-k}+\sum_k q^k \begin{bmatrix} n \\ k \end{bmatrix}\left(qs\eta^{m+1}\right)^k (x\eta)^{n+1-k}$$

$$=\sum_k \left(q^k \begin{bmatrix} n \\ k \end{bmatrix}+\begin{bmatrix} n \\ k-1 \end{bmatrix}\right)\left(qs\eta^{m+1}\right)^k (x\eta)^{n+1-k}=\sum_k \begin{bmatrix} n+1 \\ k \end{bmatrix}\left(qs\eta^{m+1}\right)^k (x\eta)^{n+1-k}.$$

Therefore we get from (2.12)

$$\begin{aligned}\varphi_m\left(F_n^*(a,b)\right)&=\sum_{k=0}^{n-1}\varphi\left(C_k^{n-1-k}(a,b)\right)=\sum_{k=0}^{n-1}\begin{bmatrix} n-k-1 \\ k \end{bmatrix}\left(qs\eta^{m+1}\right)^k (x\eta)^{n-1-2k} \\ &=\sum_{k=0}^{n-1}\begin{bmatrix} n-k-1 \\ k \end{bmatrix} q^{(m+1)\binom{k}{2}+k} s^k x^{n-1-2k}\eta^{n-1+k(m-1)}=\Phi_n(x,s,m,q)\eta^{n-1+k(m-1)}\end{aligned} \tag{5.22}$$

which implies

$$\varphi_m\left(F_n^*(a,b)\right)1=\Phi_n(x,s,m,q). \tag{5.23}$$

For the corresponding $q-$Lucas polynomials (3.18) implies the recurrence

$$\begin{aligned}&\Lambda_n(x,s,m,q)=(1+q)x\Lambda_{n-1}(x,s,m,q)+q^{n-2}(1+q)s\Lambda_{n-2}(x,q^{m-1}s,m,q) \\ &-qx^2\Lambda_{n-2}(x,s,m,q)-q^{n-2}(1+q)xs\Lambda_{n-3}(x,q^{m-1}s,m,q)-q^{2n+m-6}s^2\Lambda_{n-4}(x,q^{2m-2}s,m,q).\end{aligned} \tag{5.24}$$

A generalization of (4.24) is

$$\sum_{k=0}^{\left\lfloor \frac{n}{2}\right\rfloor}(-s)^k q^{m\binom{k}{2}}\begin{bmatrix} n \\ k \end{bmatrix}\Lambda_{n-2k}^*(x,q^{mk}s,m,q)=x^n. \tag{5.25}$$

An interesting question concerns recurrences of $\Phi_n(x,s,m,q)$ with fixed numbers $x$ and $s$. Some examples are

$$\Phi_n(x,s,-1,q) = x\Phi_{n-1}(x,s,-1,q) + q(q+1)s\Phi_{n-2}(x,s,-1,q) - q(q+1)sx\Phi_{n-3}(x,s,-1,q)$$
$$+q^{n-2}sx^2\Phi_{n-4}(x,s,-1,q) - q^3s^2\Phi_{n-4}(x,s,-1,q) + q^3s^2x\Phi_{n-4}(x,s,-1,q)$$

$$\Phi_n(x,s,0,q) = x\Phi_{n-1}(x,s,0,q) + q^{n-2}sx\Phi_{n-3}(x,s,0,q) + q^{n-2}s^2\Phi_{n-4}(x,s,0,q)$$

$$\Phi_n(x,s,1,q) = x\Phi_{n-1}(x,s,1,q) + qs\Phi_{n-2}(x,s,1,q)$$

$$\Phi_n(x,s,2,q) = [3]x\Phi_{n-1}(x,s,2,q) - q[3]x^2\Phi_{n-2}(x,s,2,q) + q^3x^3\Phi_{n-3}(x,s,2,q) - q^{2n-5}sx\Phi_{n-3}(x,s,2,q)$$
$$+q^{2n-5}sx^2\Phi_{n-4}(x,s,2,q) + q^{3n-10}s^2\Phi_{n-4}(x,s,2,q)$$

But I could not find a general rule.

## Remarks

It should be noted that analogous results hold for $q-$analogues of the binomial theorem.

The operator $U_m$ transforms the **Rogers-Szegö polynomials** $R_n(x,s) = \sum_{k=0}^{n} \begin{bmatrix} n \\ k \end{bmatrix} x^k s^{n-k}$ into

$$R_n(x,s,m) = U_m R_n(x,s) = \sum_{k=0}^{n} q^{m\binom{n-k}{2}} \begin{bmatrix} n \\ k \end{bmatrix} x^k s^{n-k}. \tag{5.26}$$

It is well known and easily verified that

$$R_n(x,s,1) = \sum_{k=0}^{n} q^{\binom{k}{2}} \begin{bmatrix} n \\ k \end{bmatrix} s^k x^{n-k} = (x+s)(x+qs)\cdots(x+q^{n-1}s) = (x+s)R_{n-1}(x,qs,1). \tag{5.27}$$

This implies that

$$R_n(x,s,m) = xR_{n-1}(x,qs,m) + sR_{n-1}(x,q^m s,m). \tag{5.28}$$

On the other hand $R_n(x,s)$ satisfies $R_{n+1}(x,s)=(x+sE)R_n(x,s),$ where $Ef(x)=f(qx).$ For

$$(x+sE)\sum_k \begin{bmatrix} n \\ k \end{bmatrix} s^k x^{n-k} = \sum_k \begin{bmatrix} n \\ k \end{bmatrix} s^k x^{n+1-k} + \sum_k \begin{bmatrix} n \\ k \end{bmatrix} s^{k+1} q^{n-k} x^{n-k}$$

$$= \sum_k \left( \begin{bmatrix} n \\ k \end{bmatrix} + q^{n-k+1} \begin{bmatrix} n \\ k-1 \end{bmatrix} \right) s^k x^{n+1-k} = \sum_k \begin{bmatrix} n+1 \\ k \end{bmatrix} s^k x^{n+1-k}.$$

Therefore we have

$$R_n(x,s,m) = xR_{n-1}(x,s,m) + sER_{n-1}(x,q^m s,m). \tag{5.29}$$

If we recall that $E=1+(q-1)xD$ we see that $R_n(x,s,m)$ satisfies the $D-$recurrence

$$R_n(x,s,m) = xR_{n-1}(x,s,m) + sR_{n-1}(x,q^m s,m) + (q-1)xsDR_{n-1}(x,q^m s,m). \tag{5.30}$$

Since $DR_n(x,s,m)=[n]R_{n-1}(x,s,m)$ this reduces to

$$R_n(x,s,m) = xR_{n-1}(x,s,m) + sR_{n-1}(x,q^m s,m) + (q^{n-1}-1)xsR_{n-2}(x,q^m s,m).$$

For $m=0$ this is the well-known recurrence for the Rogers-Szegö polynomials.